\documentclass{amsart}
\usepackage{dst}

\usepackage{ast2}

\usepackage[margin=1.4in]{geometry}

\title{A note on arithmetic in finite types}
\author{Benno van den Berg$^1$}
\address{${}^1$ ILLC, Universiteit van Amsterdam, P.O. Box 94242, 1090 GE Amsterdam, the Netherlands. E-mail: bennovdberg@gmail.com.}
\date{\today}

\begin{document}

\begin{abstract}
We present a version of arithmetic in all finite types which allows for a definition of equality at higher types for which all congruence are derivable, for which the soundness of the Dialectica interpretation is provable inside the system itself, which allows for both intensional and extensional models and for which the deduction theorem holds.
\end{abstract}

\maketitle

\section{Introduction}

Arithmetic in all finite types, or finite-type arithmetic, is a system which dates back to the work by Kreisel from the late fifties \cite{kreisel59} and has always been important in the study of constructivism. Currently, it is also playing an essential r\^ole in program extraction from from proofs and proof mining, as can be seen from the recent books \cite{schwichtenbergwainer12, kohlenbach08}. G\"odel's Dialectica interpretation is a crucial tool here. Finite-type arithmetic is also starting to attract attention in the Reverse Mathematics community, as can be seen from some recent papers on higher-order reverse mathematics like \cite{kohlenbach05c,hunter08,schweber13}.

Various versions of finite-type arithmetic exist and the differences tend to be subtle; the variety is mainly due to the fact that it is hard to find a system which has all the properties which one would like it to have. Indeed, at present the literature creates the impression that it is impossible to combine the following desirable features:
\begin{enumerate}
\item The system allows for both intensional and extensional models, such as HRO (the hereditarily recursive operations) and HEO (the hereditarily extensional operations).
\item The deduction theorem holds for this system.
\item The Dialectica interpretation is sound as an interpretation from this system into itself.
\item The system has a notion of equality at higher types which can be defined internally to the system and for which all the congruence laws are derivable.
\end{enumerate}
For example, the systems $\nha$ from \cite{troelstra73} and $\ha$ from \cite[pages 444-449]{troelstravandalen88b} have a primitive notion of equality at all finite types: for this reason atomic formulas are not decidable and this blocks the soundness of the Dialectica interpretation. The extensionality axiom can be used to reduce equality at higher types to equality at base type, which is decidable; however, the existence of a functional witnessing the Dialectica interpretation of the extensionality axiom cannot be shown inside $\eha$ itself, so this system still does not satisfy (3). In addition, it does not allow for intensional models like HRO. The intensional variant $\iha$ does satisfy (3), but it excludes models like HEO. Finally, the system $\weha$, which plays a crucial role in \cite{kohlenbach08}, excludes intensional models and works with a notion of equality for which not all congruence laws are derivable: one congruence law is valid as a rule only. Consequently, the deduction theorem fails for this system as well.

Nevertheless, it is possible to combine (1)-(4); indeed, it is the purpose of this note to introduce a version of $\ha$ which has all these desirable properties. After we have introduced it, we will see that it can be shown to be equivalent to the system called $\ha$ on page 46 of \cite{troelstra73}. The corollary that this system satisfies property (4) seems to be new.

\section{Some old versions of arithmetic in finite types}

To start, let us introduce the system called $\nha$ in \cite{troelstra73}; we will work with a formulation which includes product types (so this is the same as $\ha$ from \cite[pages 444-449]{troelstravandalen88b}).

$\nha$ is a system formulated in many-sorted intuitionistic logic, where the sorts are the finite types.
\begin{definition} The \emph{finite types} are defined by induction as follows: 0 is a finite type, and if $\sigma$ and $\tau$ are finite types, then so are $\sigma \to \tau$ and $\sigma \times \tau$. The type 0 is the \emph{ground} or \emph{base type}, while the other types will be called \emph{higher types}.
\end{definition}
There will be infinitely many variables of each sort. In addition, there will be constants:
\begin{enumerate}
\item for each pair of types $\sigma, \tau$ a combinator ${\bf k}^{\sigma, \tau}$ of sort $\sigma \to (\tau \to \sigma)$.
\item for each triple of types $\rho, \sigma, \tau$ a combinator ${\bf s}^{\rho, \sigma, \tau}$ of type $(\rho \to (\sigma \to \tau)) \to ((\rho \to \sigma) \to (\rho \to \tau))$.
\item for each pair of types $\rho, \sigma$ combinators ${\bf p}^{\rho, \sigma}, {\bf p}^{\rho, \sigma}_0, {\bf p}^{\rho, \sigma}_1$ of types $\rho \to (\sigma \to \rho \times \sigma)$, $\rho \times \sigma \to \rho$ and $\rho \times \sigma \to \sigma$, respectively.
\item a constant 0 of type 0 and a constant $S$ of type $0 \to 0$.
\item for each type $\sigma$ a combinator ${\bf R}^\sigma$ (``the recursor'') of type $\sigma \to ((0 \to (\sigma \to \sigma)) \to (0 \to \sigma))$.
\end{enumerate}

\begin{definition}
The terms of $\nha$ are defined inductively as follows:
\begin{itemize}
\item each variable or constant of type $\sigma$ will be a term of type $\sigma$.
\item if $f$ is a term of type $\sigma \to \tau$ and $x$ is a term of type $\sigma$, then $fx$ is a term of type $\tau$.
\end{itemize}
\end{definition}
The convention is that application associates to the left, which means that an expression like $fxyz$ has to be read as $(((fx)y)z)$.

\begin{definition}
The formulas are defined inductively as follows:
\begin{itemize}
\item $\bot$ is a formula and if $s$ and $t$ are terms of the same type $\sigma$, then $s =_\sigma t$ is a formula.
\item if $\varphi$ and $\psi$ are formulas, then so are $\varphi \land \psi, \varphi \lor \psi, \varphi \to \psi$.
\item if $x$ is a variable of type $\sigma$ and $\varphi$ is a formula, then $\exists x^\sigma \, \varphi$ and $\forall x^\sigma \, \varphi$ are formulas.
\end{itemize}
\end{definition}

Finally, the axioms and rules of $\nha$ are:
\begin{enumerate}
\item[(i)] All the axioms and rules of many-sorted intuitionistic logic (say in Hilbert-style).
\item[(ii)] Equality is an equivalence relation at all types:
\[ x = x, \qquad x = y \to y = x, \qquad x = y \land y = z \to x = z \]
\item[(iii)] The congruence laws for equality at all types:
\[ f = g \to fx = gx, \qquad x = y \to fx = fy \]
\item[(v)] The successor axioms:
\[ \lnot S(x) = 0, \qquad S(x) = S(y) \rightarrow x = y \]
\item[(v)] For any formula $\varphi$ in the language of $\ha$, the induction axiom:
    \[ \varphi(0) \to \big( \, \forall x^0 \, ( \, \varphi(x) \to \varphi(Sx) \, ) \to \forall x^0 \, \varphi(x) \, \big). \]
\item[(vi)] The axioms for the combinators:
\begin{eqnarray*}
{\bf k} x y & = & x \\
{\bf s} x y z & = & xz(yz) \\
{\bf p}_0({\bf p}xy) & = & x \\
{\bf p}_1({\bf p}xy) & = & y \\
{\bf p}({\bf p}_0x)({\bf p}_1x) & = & x
\end{eqnarray*} as well as for the recursor:
\begin{eqnarray*}
{\bf R} x y 0 & = & x \\
{\bf R} x y (Sn) & = & yn({\bf R}xyn)
\end{eqnarray*}
\end{enumerate}

The system $\eha$ is obtained from $\nha$ by adding the axiom of extensionality:
\[ \EXT: \qquad \forall f^{\sigma \to \tau}, g^{\sigma \to \tau} \big( \, (\forall x^{\sigma} fx =_\tau gx) \to f =_{\sigma \to \tau} g \big ); \]
on the other hand, the system $\iha$ adds to $\nha$ combinators ${\bf e}^\sigma$ of type $\sigma \to (\sigma \to 0)$ and axioms
\[ {\bf e}^\sigma xy \leq 1, \qquad {\bf e}^\sigma xy = 0 \leftrightarrow x =_\sigma y. \]

\section{A new version of arithmetic in finite types}

As mentioned in the introduction, the system $\nha$ has the problem that it has a primitive notion of equality at higher types: consequently, atomic formulas are not decidable and this blocks the soundness of the Dialectica interpretation. We will try to solve this by defining equality at higher types in terms of equality at ground type. This is presumably not possible in $\nha$, but it is possible in both $\eha$ and $\iha$: in $\eha$ one can define equality at higher types extensionally, while in $\iha$ one can directly use the combinator $\bf e$ to reduce equalities at higher types to equalities at base type. Both these systems, however, exclude important classes of models (like HRO in case of $\eha$ and HEO in case of $\iha$) and the way they reduce equalities at higher types to equalities at base type are often incompatible. Therefore we are looking for a ``modular'' way to reduce equalities at higher type to equalities at base type which is compatible with both extensional and intensional features.

To see how to do this, let us consider the following \emph{principle of observational equivalence}:
\[ \OBS \qquad \forall f^{\sigma \to 0} ( \, fx =_0 fy \, ) \to x =_\sigma y. \]
\begin{proposition}
We have $\eha \vdash \OBS$ and $\iha \vdash \OBS$, while $\nha \not\vdash \OBS$.
\end{proposition}
\begin{proof}
Let $x$ and $y$ be two objects of type $\sigma$ such that $\forall f^{\sigma \to 0} \, fx = fy$. First we work in $\eha$. Without loss of generality, we may assume that $\sigma$ is of the form $\rho \to 0$. But then we have for any $z$ of type $\rho$ that
\[ xz = (\lambda a.az) x = (\lambda a.az) y = yz, \]
using our assumption with $f^{\sigma \to 0} = \lambda a.az$; so $x = y$ by the extensionality axiom.

Next, we work in $\iha$. In this case we simply take $f = {\bf e}^\sigma x$. Then $fx = 0$, so if $fx = fy$, then $fy =0$ and $x = y$.

If, on the other hand, we take HEO as our model of G\"odel's $\T$, but we take equality of codes as our notion of equality (as in HRO), then we get a model of $\nha$ (even $\npa$) in which $\OBS$ fails: for in this model two different codes of the same total recursive function will be observationally equivalent as type 1 objects, but will not be equal in the sense of the model.
\end{proof}

So the principle of observational equivalence reduces equality at higher types to equality at base type in a way which is valid both on the extensional and the intensional viewpoint. From a philosophical point of view one may argue in favour of this principle as follows: although higher-type objects can be regarded as finitary in some sense, they remain fairly abstract, as opposed to objects of type 0, which are concrete natural numbers. Therefore it makes sense to say for objects of higher type that they should be regarded as the same if under any attempt to extract something concrete, i.e., a natural number, from them, they yield the same result (in that sense, they are ``observationally equivalent'').

Let us now define our new version of Heyting arithmetic in all finite types, which we will call $\hazero$. $\hazero$ is a system formulated in many-sorted intuitionistic logic, where the sorts are the finite types. There will be infinitely many variables of each sort. In addition, there will be constants:
\begin{enumerate}
\item for each pair of types $\sigma, \tau$ a combinator ${\bf k}^{\sigma, \tau}$ of sort $\sigma \to (\tau \to \sigma)$.
\item for each triple of types $\rho, \sigma, \tau$ a combinator ${\bf s}^{\rho, \sigma, \tau}$ of type $(\rho \to (\sigma \to \tau)) \to ((\rho \to \sigma) \to (\rho \to \tau))$, as well as a combinator ${\bf b}^{\rho,\sigma,\tau}$ of type $(\sigma \to \tau) \to ((\rho \to \sigma) \to (\rho \to \tau))$, as well as a combinator ${\bf q}^{\rho,\sigma, \tau}$ of type $(\sigma \to \tau) \to (\rho \to ((\rho \to \sigma) \to \tau))$.
\item for each pair of types $\rho, \sigma$ combinators ${\bf p}^{\rho, \sigma}, {\bf p}^{\rho, \sigma}_0, {\bf p}^{\rho, \sigma}_1$ of types $\rho \to (\sigma \to \rho \times \sigma)$, $\rho \times \sigma \to \rho$ and $\rho \times \sigma \to \sigma$, respectively.
\item a constant 0 of type 0 and a constant $S$ of type $0 \to 0$.
\item for each type $\sigma$ a combinator ${\bf R}^\sigma$ (``the recursor'') of type $\sigma \to ((0 \to (\sigma \to \sigma)) \to (0 \to \sigma))$.
\end{enumerate}
Note that we have added two additional combinators: $\bf b$ and $\bf q$. The reason should become clear shortly (see Remark \ref{oversight} below).
\begin{definition}
The terms of $\hazero$ are defined inductively as follows:
\begin{itemize}
\item each variable or constant of type $\sigma$ will be a term of type $\sigma$.
\item if $f$ is a term of type $\sigma \to \tau$ and $x$ is a term of type $\sigma$, then $fx$ is a term of type $\tau$.
\end{itemize}
\end{definition}
The convention is that application associates to the left, which means that an expression like $fxyz$ has to be read as $(((fx)y)z)$.

\begin{definition}
The formulas of $\hazero$ are defined inductively as follows:
\begin{itemize}
\item $\bot$ is a formula and if $s$ and $t$ are terms of the type 0, then $s =_0 t$ is a formula.
\item if $\varphi$ and $\psi$ are formulas, then so are $\varphi \land \psi, \varphi \lor \psi, \varphi \to \psi$.
\item if $x$ is a variable of type $\sigma$ and $\varphi$ is a formula, then $\exists x^\sigma \, \varphi$ and $\forall x^\sigma \, \varphi$ are formulas.
\end{itemize}
\end{definition}

Equality at higher types will be defined ``observationally'', as follows:
\begin{eqnarray*}
x =_{\sigma} y & := & \forall f^{\sigma \to 0} \, fx =_0 fy.
\end{eqnarray*}

Finally, the axioms and rules of $\hazero$ are:
\begin{enumerate}
\item[(i)] All the axioms and rules of many-sorted intuitionistic logic (say in Hilbert-style).
\item[(ii)] Equality at type 0 is an equivalence relation:
\[ x =_0 x, \qquad x =_0 y \to y =_0 x, \qquad x =_0 y \land y =_0 z \to x =_0 z \]
\item[(iii)] There is one additional congruence law:
\[ x =_0 y \to fx =_0 fy \]
\item[(v)] There are successor axioms:
\[ \lnot S(x) =_0 0, \qquad S(x) =_0 S(y) \rightarrow x =_0 y \]
\item[(v)] For any formula $\varphi$ in the language of $\hazero$, the induction axiom:
    \[ \varphi(0) \to \big( \, \forall x^0 \, ( \, \varphi(x) \to \varphi(Sx) \, ) \to \forall x^0 \, \varphi(x) \, \big). \]
\item[(vi)] The axioms for the combinators
\begin{eqnarray*}
{\bf k} x y & = & x \\
{\bf s} x y z & = & xz(yz) \\
{\bf b}xyz & = & x(yz) \\
{\bf q}xyz & = & x(zy) \\
{\bf p_0}({\bf p}xy) & = & x \\
{\bf p_1}({\bf p}xy) & = & y \\
{\bf p}({\bf p}_0x)({\bf p}_1x) & = & x
\end{eqnarray*} as well as for the recursor:
\begin{eqnarray*}
{\bf R} x y 0 & = & x \\
{\bf R} x y (Sn) & = & yn({\bf R}xyn);
\end{eqnarray*}
here equality means observational equivalence, as defined above.
\end{enumerate}

We have defined equality at higher types as observational equivalence; but calling observational equivalence equality does not make it act like equality. Therefore the first thing we need to do is to prove that in $\hazero$ observational equivalence acts as a congruence. To this purpose, note that we can define combinators ${\bf i} := {\bf skk}$ and ${\bf t} := {\bf qi}$, for which we can derive
\[ \begin{array}{l}
{\bf i}x = {\bf skk}x = {\bf k}x({\bf k}x) = x, \\
{\bf t}xy = {\bf qi}xy = {\bf i}(yx) = yx,
\end{array} \]
without using any congruence laws.

\begin{proposition} $\hazero \vdash x =_\sigma y \to f x =_\tau fy$ and $\hazero \vdash f =_{\sigma \to \tau} g \to fx =_\tau gx$.
\end{proposition}
\begin{proof}
Note that we have $\hazero \vdash x =_\sigma y \to f x =_0 fy$: if $\sigma = 0$, then this is an axiom; if $\sigma$ is a higher type, it holds by definition. So we only need to prove $\hazero \vdash x =_\sigma y \to f x =_\tau fy$ in case $\tau$ is a higher type. In that case, let $u^{\tau \to 0}$ be arbitrary and consider the term ${\bf b}uf$ of type $\sigma \to 0$. Since $x =_\sigma y$, we have
\[ u(fx) = {\bf b}ufx = {\bf b}ufy = u(fy), \]
so $fx = fy$, by definition of equality as observational equivalence.

To prove the congruence axiom $f =_{\sigma \to \tau} g \to fx =_\tau gx$ in $\hazero$, we make a case distinction:
\begin{itemize}
\item $\tau = 0$: note that ${\bf t}x$ is of type $(\sigma \to 0) \to 0$, so $f =_{\sigma \to 0} g$ implies ${\bf t}xf = {\bf t}xg$ and hence $fx = gx$.
\item If $\tau$ is of higher type, we need to show that $f =_{\sigma \to \tau} g$ implies $ufx = ugx$ for any $u$ of type $\tau \to 0$. To this purpose, consider the term ${\bf q}ux$ of type $(\sigma \to \tau) \to 0$. From $f =_{\sigma \to \tau} g$ it follows that  ${\bf q}uxf = {\bf q}uxg$ and hence $u(fx) = u(gx)$.
\end{itemize}
\end{proof}

Now that we have shown that observational equivalence is a congruence we can define $\lambda$-abstraction using $\bf k$ and $\bf s$, as in \cite[Proposition 9.1.8]{troelstravandalen88b}, for instance. Note that we can do this only now, because the proof that the $\lambda$-abstraction defined using $\bf k$ and $\bf s$ acts as it should (i.e., proves $\beta$-equality) \emph{uses the congruence laws for equality}.

\begin{remark} \label{oversight}
{\rm This is overlooked on pages 452 and 453 of \cite{troelstravandalen88b}: one cannot freely use $\lambda$-abstractions to prove the congruence laws, because the congruence laws are used in the proof of the combinatory completeness of $\bf k$ and $\bf s$. (The argument is ascribed to Rath, but, as he was working with a version of $\nha$ based on the $\lambda$-calculus rather than on combinatory logic, his argument was not circular.) The necessity of breaking this circle was our reason for introducing the additional combinators $\bf b$ and $\bf q$. We feel that the system we have called $\hazero$ achieves what the system called $\hazero$ on pages 452 and 453 of \cite{troelstravandalen88b} was meant to achieve, so it seems appropriate to use the same name.}
\end{remark}

\begin{corollary}
The following systems all prove the same theorems in their common language: $\hazero$, $\nha$, $\nha + \OBS$ and the system called $\ha$ on page 46 of \cite{troelstra73}.
\end{corollary}
\begin{proof}
The only observation to make is that once one has combinatory completeness, as one has in $\nha$, $\nha + \OBS$ and $\ha$, one can define $\bf b$ and $\bf q$. For example, one could put
\begin{eqnarray*}
{\bf b} & := & {\bf s}({\bf ks}){\bf k}, \\
{\bf q} & := & {\bf b}({\bf s}({\bf bbs})({\bf kk})){\bf b},
\end{eqnarray*}
as one may verify.
\end{proof}

This means that $\hazero$ allows for both intensional and extensional models, such as HRO and HEO; also the deduction theorem holds for $\hazero$, because, unlike $\weha$, it extends multi-sorted first-order intuitionistic logic with axioms only.

In addition, the Dialectica interpretation works as an interpretation of $\hazero$ in $\hazero$, essentially because observational equality is defined using a universal formula. This means that all axioms in groups (ii), (iii), (iv) and (vi) are universal and hence interpreted by themselves. This should be compared with what happens in $\eha$: if one defines equality at higher types extensionally, then the extensionality axiom turns into a congruence law. This congruence law, however, is not universal and its Dialectica interpretation is not trivial (in fact, as shown by Howard in the appendix of \cite{troelstra73}, the existence of a realizer cannot be shown in $\bf ZF$ set theory). Finally, also the nonstandard Dialectica interpretation from \cite{bergbriseidsafarik12} works for this system: we assumed extensionality throughout, but the only thing which is needed for the soundness proof is that the congruence laws for equality hold. For this reason the system $\hazero$ might prove useful if one wishes to combine methods from \cite{bergbriseidsafarik12} with traditional proof mining techniques.

\bibliographystyle{plain} \bibliography{dst}

\end{document}